\documentclass[11pt]{article}
\usepackage{mathrsfs}
\usepackage{amssymb}
\usepackage{amsmath}
\usepackage{latexsym}
\usepackage{amsfonts}

\usepackage[all]{xy}
\def\Coker{\mathop{\rm Coker}\nolimits}

\def\mod{\mathop{\rm mod}\nolimits}
\def\rad{\mathop{\rm rad}\nolimits}
\def\grade{\mathop{\rm grade}\nolimits}

\def\id{\mathop{\rm id}\nolimits}
\def\pd{\mathop{\rm pd}\nolimits}
\def\gldim{\mathop{\rm gl.dim}\nolimits}
\def\add{\mathop{\rm add}\nolimits}

\title{\Large \bf Higher Auslander Algebras Admitting Trivial Maximal Orthogonal Subcategories
}
\author{Zhaoyong Huang\thanks{{\it E-mail address}:
huangzy@nju.edu.cn}, Xiaojin Zhang\thanks{{\it E-mail address}:
xiaojinzhang@sohu.com}\\
{\footnotesize \it Department of Mathematics, Nanjing University,
Nanjing 210093, Jiangsu Province, P. R. China}}
\date{ }

\begin{document}

\baselineskip=18pt

\maketitle

\begin{abstract}
For an Artinian $(n-1)$-Auslander algebra $\Lambda$ with global
dimension $n(\geq 2)$, we show that if $\Lambda$ admits a trivial
maximal $(n-1)$-orthogonal subcategory of $\mod\Lambda$, then
$\Lambda$ is a Nakayama algebra and the projective or injective
dimension of any indecomposable module in $\mod\Lambda$ is at most
$n-1$. As a result, for an Artinian Auslander algebra with global
dimension 2, if $\Lambda$ admits a trivial maximal $1$-orthogonal
subcategory of $\mod\Lambda$, then $\Lambda$ is a tilted algebra of
finite representation type. Further, for a finite-dimensional
algebra $\Lambda$ over an algebraically closed field $K$, we show
that $\Lambda$ is a basic and connected $(n-1)$-Auslander algebra
$\Lambda$ with global dimension $n(\geq 2)$ admitting a trivial
maximal $(n-1)$-orthogonal subcategory of $\mod\Lambda$ if and only
if $\Lambda$ is given by the quiver:
$$\xymatrix{1 &
\ar[l]_{\beta _{1}} 2 & \ar[l]_{\beta _{2}} 3 & \ar[l]_{\beta _{3}}
\cdots & \ar[l]_{\beta _{n}} n+1}
$$ modulo the ideal generated by $\{
\beta_{i}\beta_{i+1}| 1\leq i\leq n-1 \}$. As a consequence, we get
that a finite-dimensional algebra over an algebraically closed field
$K$ is an $(n-1)$-Auslander algebra with global dimension $n(\geq
2)$ admitting a trivial maximal $(n-1)$-orthogonal subcategory if
and only if it is a finite direct product of $K$ and $\Lambda$ as
above. Moreover, we give some necessary condition for an Artinian
Auslander algebra admitting a non-trivial maximal 1-orthogonal
subcategory.

\end{abstract}

\vspace{0.5cm}

\centerline{\large \bf 1. Introduction}

\vspace{0.2cm}

It is well known that the notion of maximal $n$-orthogonal
subcategories introduced by Iyama in [Iy3] played a crucial role in
developing the higher-dimensional Auslander-Reiten theory (see [Iy3]
and [Iy4]). This notion coincides with that of $(n+1)$-cluster
tilting subcategories introduced by Keller and Reiten in [KR]. In
general, maximal $n$-orthogonal subcategories rarely exist. So it
would be interesting to investigate when maximal $n$-orthogonal
subcategories exist and the properties of algebras admitting such
subcategories. Several authors have worked on this topics (see [EH],
[GLS], [HuZ], [Iy3], [Iy4], [Iy5], [Iy6], [L], and so on). As a
generalization of the notion of the classical Auslander algebras,
Iyama introduced the notion of $n$-Auslander algebras in [Iy5]. Then
he proved that for a finite-dimensional $(n-1)$-Auslander algebra
$\Lambda$ with global dimension $n (\geq 2)$ over an algebraically
closed field $K$, $\Lambda$ has maximal $(n-1)$-orthogonal modules
in $\mod \Lambda$ if and only if $\Lambda$ is Morita equivalent to
$T^{(n)}_{m}(K)$ for some $m \geq 1$, where $T^{(1)}_{m}(K)$ is an
$m\times m$ upper triangular algebra and $T^{(n)}_{m}(K)$ is the
endomorphism algebra of a maximal $(n-2)$-orthogonal module in $\mod
T^{(n-1)}_{m}(K)$. Moreover, he gave some examples of the quivers of
these algebras inductively. In [HuZ] we proved that an Artinian
$(n-1)$-Auslander algebra $\Lambda$ with global dimension $n (\geq
2)$ admits a trivial maximal $(n-1)$-orthogonal subcategory of $\mod
\Lambda$ if and only if any simple module $S\in \mod \Lambda$ with
projective dimension $n$ is injective. In [HuZ] we also proved that
for an almost hereditary algebra $\Lambda$ with global dimension 2,
if $\Lambda$ admits a maximal $1$-orthogonal subcategory
$\mathscr{C}$ of $\mod\Lambda$, then $\mathscr{C}$ is trivial. In
this paper, we continue to study the structure of an
$(n-1)$-Auslander algebra $\Lambda$ admitting a trivial maximal
$(n-1)$-orthogonal subcategory of $\mod\Lambda$. This paper is
organized as follows.

In Section 2, we give some notions and notations and collect some
preliminary results about minimal morphisms. In Section 3, we give
some homological properties of indecomposable modules (in
particular, simple modules) over higher Auslander algebras
(admitting a trivial maximal orthogonal subcategory of $\mod
\Lambda$).

Let $\Lambda$ be an Artinian $(n-1)$-Auslander algebra $\Lambda$
with global dimension $n (\geq 2)$ admitting a trivial maximal
$(n-1)$-orthogonal subcategory of $\mod\Lambda$. We get in Section 4
the following conclusions: (1) For any simple module $S\in
\mod\Lambda$ with $\pd_{\Lambda}S=n$, the $i$th syzygy module of $S$
is simple $1\leq i\leq n$ and all terms in a minimal projective
resolution of $S$ are indecomposable. (2) The sum of the projective
and injective dimensions of any non-projective-injective
indecomposable module in $\mod \Lambda$ is equal to $n$. As a
consequence, the projective or injective dimension of any
indecomposable module in $\mod\Lambda$ is at most $n-1$. (3)
$\Lambda$ is a Nakayama algebra. As a result, we get that an
Artinian Auslander algebra $\Lambda$ admitting a trivial maximal
1-orthogonal subcategory of $\mod \Lambda$ is a tilted algebra of
finite representation type.

As an application of the results obtained in Section 4, we give in
Section 5 an explicit description of the quiver of a
finite-dimensional $(n-1)$-Auslander algebra $\Lambda$ with global
dimension $n$ admitting a trivial maximal $(n-1)$-orthogonal
subcategory of $\mod \Lambda$. Let $K$ be an algebraically closed
field. We first prove that for a basic and connected
finite-dimensional $(n-1)$-Auslander $K$-algebra $\Lambda$, if
$\Lambda$ admits a trivial maximal $(n-1)$-orthogonal subcategory of
$\mod\Lambda$, then there exists a unique simple module
$S\in\mod\Lambda$ with projective dimension $n$. Then we prove that
$\Lambda$ is a basic and connected finite-dimensional
$(n-1)$-Auslander $K$-algebra admitting a trivial maximal
$(n-1)$-orthogonal subcategory of $\mod\Lambda$ if and only if
$\Lambda$ is given by the quiver:
$$\xymatrix{1 & \ar[l]_{\beta _{1}} 2 & \ar[l]_{\beta _{2}} 3 &
\ar[l]_{\beta _{3}} \cdots & \ar[l]_{\beta _{n}} n+1}
$$ modulo the ideal generated by $\{
\beta_{i}\beta_{i+1}| 1\leq i\leq n-1 \}$. As a consequence, we
establish the following structure theorem: a finite-dimensional
$(n-1)$-Auslander $K$-algebra with global dimension $n (\geq 2)$
admitting a trivial maximal $(n-1)$-orthogonal subcategory of $\mod
\Lambda$ if and only if it is a finite direct product of $K$ and
$\Lambda$ as above.

By [Iy5], there exist an Auslander algebra with global dimension 2
admitting a non-trivial maximal 1-orthogonal subcategory. On the
other hand, by [HuZ, Corollary 3.12] we have that if $\Lambda$ is an
Artinian Auslander algebra with global dimension 2 admitting a
non-trivial maximal 1-orthogonal subcategory of $\mod\Lambda$, then
there exists a simple module $S\in\mod \Lambda$ such that both the
projective and injective dimensions of $S$ are equal to 2. In
Section 6, we further give a necessary condition for an Auslander
algebra with global dimension 2 admitting a non-trivial maximal
1-orthogonal subcategory in terms of the homological properties of
simple modules. We prove that if $\Lambda$ is an Artinian Auslander
algebra with global dimension 2 admitting a non-trivial maximal
1-orthogonal subcategory of $\mod\Lambda$, then there exist at least
two non-injective simple modules in $\mod \Lambda$ with projective
dimension 2.

\vspace{0.5cm}

\centerline{\large \bf 2. The properties of minimal morphisms}

\vspace{0.2cm}

In this section, we give some notions and notations in our
terminology and collect some preliminary results about minimal
morphisms for later use.

Throughout this paper, $\Lambda$ is an Artinian algebra with the
center $R$, $\mod \Lambda$ is the category of finitely generated
left $\Lambda$-modules and $\gldim \Lambda$ denotes the global
dimension of $\Lambda$. We denote by $\mathbb{D}$ the ordinary
duality, that is, $\mathbb{D}(-)={\rm Hom}_R(-, I(R/J(R)))$, where
$J(R)$ is the Jacobson radical of $R$ and $I(R/J(R))$ is the
injective envelope of $R/J(R)$.

Let $M$ be in $\mod \Lambda$. We use
$$\cdots \rightarrow P_{i}(M)\rightarrow \cdots \rightarrow
P_{1}(M)\rightarrow P_{0}(M)\rightarrow M \rightarrow 0$$ and
$$0 \rightarrow M\rightarrow I^{0}(M)\rightarrow
I^{1}(M)\rightarrow \cdots \rightarrow I^{i}(M) \rightarrow \cdots$$
to denote the minimal projective resolution and the minimal
injective resolution of $M$, respectively. In particular, $P_{0}(M)$
and $I^{0}(M)$ are the projective cover and the injective envelope
of $M$, respectively. Denote by $\Omega^{i}M$ and $\Omega^{-i}M$ the
$i$th syzygy and co-syzygy of $M$, respectively.

The following easy observations are well-known.

\vspace{0.2cm}

{\bf Lemma 2.1} {\it Let $M\in \mod \Lambda$ and $M\cong
M_{1}\bigoplus M_{2}$. Then $$0\rightarrow M(\cong M_{1}\bigoplus
M_{2})\rightarrow I^{0}(M^{'})\bigoplus I^{0}(M^{''})\rightarrow
I^{1}(M^{'})\bigoplus I^{1}(M^{''})\rightarrow\cdots$$ and
$$\cdots\rightarrow P_{1}(M^{'})\bigoplus P_{1}(M^{''})\rightarrow
P_{0}(M^{'})\bigoplus P_{0}(M^{''})\rightarrow M(\cong
M_{1}\bigoplus M_{2})\rightarrow 0$$ are a minimal injective
resolution and a minimal projective resolution of $M$, respectively,
and $\Omega^{-i}M\cong \Omega^{-i}M_{1}\bigoplus\Omega^{-i}M_{2}$
and $\Omega^{i}M\cong \Omega^{i}M_{1}\bigoplus\Omega^{i}M_{2}$ for
any $i\geq 1$.}

\vspace{0.2cm}

{\bf Lemma 2.2} {\it Let $M$ and $S$ be in $\mod \Lambda$ with $S$
simple. Then ${\rm Ext}_{\Lambda}^{i}(S, M)\cong {\rm
Hom}_{\Lambda}(S,\Omega^{-i}M)$ for any $i \geq 0$.}

\vspace{0.2cm}

Recall from [AuR] that a morphism $f: M\rightarrow N$ in $\mod
\Lambda$ is said to be {\it left minimal} if an endomorphism
$g:N\rightarrow N$ is an automorphism whenever $f=gf$. Dually, the
notion of right minimal morphisms is defined.

\vspace{0.2cm}

{\bf Lemma 2.3} ([Au, Chapter II, Lemma 4.3]) {\it Let $0\rightarrow
A\stackrel{g}{\rightarrow }B\stackrel{f}{\rightarrow} C\rightarrow
0$ be a non-split exact sequence in $\mod\Lambda$.

(1) If $A$ is indecomposable, then $f:B\rightarrow C$ is right
minimal.

(2) If $C$ is indecomposable, then $g:A\rightarrow B$ is left
minimal.}

\vspace{0.2cm}

By Lemma 2.3, we immediately have the following result.

\vspace{0.2cm}

{\bf Corollary 2.4} {\it Let $M\in \mod \Lambda$ be an
indecomposable non-injective module and $I^{0}(M)$ projective. Then
$$\cdots \rightarrow P_{i}(M)\rightarrow \cdots \rightarrow
P_{1}(M)\rightarrow P_{0}(M)\rightarrow I^{0}(M)
\stackrel{\pi}{\rightarrow} I^{0}(M)/M \rightarrow 0$$ is a minimal
projective resolution of $I^{0}(M)/M$, where $\pi$ is the natural
epimorphism.}

\vspace{0.2cm}

The following properties of minimal morphisms are useful in the rest
of the paper.

\vspace{0.2cm}

{\bf Lemma 2.5} {\it Let $0\rightarrow
A\stackrel{g}{\rightarrow}B\stackrel{f}{\rightarrow} C\rightarrow 0$
be a non-split exact sequence in $\mod\Lambda$.

(1) If $g$ is left minimal, then ${\rm Ext}_{\Lambda}^{1}(C^{'},
A)\neq 0$ for any non-zero direct summand $C^{'}$ of $C$.

(2) If $f$ is right minimal, then ${\rm Ext}_{\Lambda}^{1}(C,
A^{'})\neq 0$ for any non-zero direct summand $A^{'}$ of $A$.}

\vspace{0.2cm}

{\it Proof.} (1) If ${\rm Ext}_{\Lambda}^{1}(C^{'}, A)=0$ holds for
some non-zero direct summand $C^{'}$ of $C$. Then we have the
following commutative diagram:

$$\xymatrix{0 \ar[r] & A \ar@{=}[d] \ar[r] & C^{'}\bigoplus A
\ar[d]^{i_1}
\ar[r]^{\pi _3} & C^{'} \ar[d]^{i_2} \ar[r] \ar@{-->}@<2pt>[l]^{i_3}& 0\\
0 \ar[r] & A \ar[r]^{g} & B \ar[r]^{f} & C \ar[r]
\ar@{-->}@<2pt>[u]^{\pi _2}& 0}
$$ such that $\pi_{3}i_{3}=1_{C^{'}}=\pi_{2} i_{2}$
and $i_{2}\pi_{3}=f i_{1}$. Then $1_{C^{'}}=(\pi_{2}
i_{2})(\pi_{3}i_{3})=(\pi_{2}f)(i_{1}i_{3})$, and hence $C^{'}$ is a
direct summand of $B$ and $(\pi_{2}f)g=0$. By [AuRS, Chapter I,
Theorem 2.4], $g$ is not left minimal, which is a contradiction.

Similarly, we get (2). $\hfill{\square}$

\vspace{0.2cm}

The following lemma establishes a connection between left minimal
morphisms and right minimal morphisms.

\vspace{0.2cm}

{\bf Lemma 2.6} {\it Let $$0\rightarrow A\stackrel{g}{\rightarrow
}B\stackrel{f}{\rightarrow} C\rightarrow 0 \eqno{(1)}$$ be a
non-split exact sequence in $\mod\Lambda$ with $B$
projective-injective. Then the following are equivalent.

(1) $A$ is indecomposable and $g$ is left minimal.

(2) $C$ is indecomposable and $f$ is right minimal.}

\vspace{0.2cm}

{\it Proof.} $(1)\Rightarrow(2)$ Since $A$ is indecomposable, $f$ is
right minimal by Lemma 2.3. Then $B$ is projective implies that the
exact sequence (1) is part of a minimal projective resolution of
$C$. If $C=C_{1}\bigoplus C_{2}$ with $C_1$ and $C_2$ non-zero, then
neither $C_{1}$ nor $C_{2}$ are projective by Lemma 2.5. So both
$\Omega^{1}C_{1}$ and $\Omega^{1}C_{2}$ are non-zero and $A\cong
\Omega^{1}C_{1}\bigoplus \Omega^{1}C_{2}$, which contradicts with
that $A$ is indecomposable.

Similarly, we get $(2)\Rightarrow (1)$. $\hfill{\square}$

\vspace{0.5cm}

\centerline{\large \bf 3. Higher Auslander algebras and maximal
orthogonal subcategories}

\vspace{0.2cm}

In this section, we give the definitions of higher Auslander
algebras and maximal orthogonal subcategories, which were introduced
by Iyama in [Iy5] and [Iy3], respectively. Then we study the
homological behavior of indecomposable modules (in particular,
simple modules) over higher Auslander algebras (admitting a trivial
maximal orthogonal subcategory of $\mod \Lambda$).

As a generalization of the notion of classical Auslander algebras,
Iyama introduced in [Iy5] the notion of $n$-Auslander algebras as
follows.

\vspace{0.2cm}

{\bf Definition 3.1} ([Iy5]) For a positive integer $n$, $\Lambda$
is called an {\it $n$-Auslander algebra} if $\gldim \Lambda\leq n+1$
and $I^{0}(\Lambda),I^{1}(\Lambda), \cdots, I^{n}(\Lambda)$ are
projective.

\vspace{0.2cm}

The notion of $n$-Auslander algebras is left and right symmetric by
[Iy5, Theorem 1.10]. It is trivial that $n$-Auslander algebras with
global dimension at most $n$ are semisimple. In particular, the
notion of 1-Auslander algebras is just that of classical Auslander
algebras. In the following, we assume that $n \geq 2$ when an
$(n-1)$-Auslander algebra is concerned.

 Denote by
$\mathscr{PI}^{n}(\Lambda)$ (resp. $\mathscr{IP}^{n}(\Lambda)$) the
subcategory of $\mod \Lambda$ consisting of indecomposable
projective modules with injective dimension $n$ (resp.
indecomposable injective modules with projective dimension $n$). By
applying Lemma 2.6 to $(n-1)$-Auslander algebras, we get the
following result.

\vspace{0.2cm}

{\bf Lemma 3.2} {\it Let $\Lambda$ be an $(n-1)$-Auslander algebra
with $\gldim\Lambda=n$. Then we have the following

(1) For any $P \in \mathscr{PI}^{n}(\Lambda)$, the minimal injective
resolution of $P$:
$$0\rightarrow P\rightarrow I^{0}(P)\rightarrow
I^{1}(P)\rightarrow\cdots \rightarrow I^{n}(P)\rightarrow 0
\eqno{(2)}$$ is a minimal projective resolution of $I^{n}(P)$ and
$I^{n}(P)$ is indecomposable.

(2) For any module $I \in \mathscr{IP}^{n}(\Lambda)$, the minimal
projective resolution of $I$:
$$0\rightarrow P_{n}(I)\rightarrow \cdots \rightarrow P_{1}(I)
\rightarrow P_{0}(I) \rightarrow I \rightarrow 0$$ is a minimal
injective resolution of $P_{n}(I)$ and $P_{n}(I)$ is
indecomposable.}

\vspace{0.2cm}

{\it Proof.} (1) Since $\Lambda$ is an $(n-1)$-Auslander algebra, by
Lemma 2.1, it is easy to see that $I^{i}(P)$ is projective for any
$0\leq i\leq n-1$. So the exact sequence $(2)$ is a projective
resolution of $I^{n}(P)$, and then the assertion follows from Lemma
2.6.

Dually, we get (2). $\hfill{\square}$

\vspace{0.2cm}

By Lemma 3.2, we get immediately the following result.

\vspace{0.2cm}

{\bf Lemma 3.3} {\it Let $\Lambda$ be an $(n-1)$-Auslander algebra
with $\gldim\Lambda=n$. Then $\Omega^{n}$ gives a one-one
correspondence between $\mathscr{IP}^{n}(\Lambda)$ and
$\mathscr{PI}^{n}(\Lambda)$ with the inverse $\Omega^{-n}$.}

\vspace{0.2cm}

For a module $M\in \mod \Lambda$, we use $\pd_{\Lambda}M$ and $\id
_{\Lambda}M$ to denote the projective dimension and the injective
dimension of $M$, respectively.

\vspace{0.2cm}

{\bf Lemma 3.4} {\it Let $\Lambda$ be an $(n-1)$-Auslander algebra
with $\gldim\Lambda=n$ and $S\in \mod\Lambda$ a simple module with
$\pd_{\Lambda}S=n$, then $P_{n}(S)$ is indecomposable.}

\vspace{0.2cm}

{\it Proof.} Let $\Lambda$ be an $(n-1)$-Auslander algebra with
$\gldim\Lambda=n$ and $S\in \mod\Lambda$ a simple module with
$\pd_{\Lambda}S=n$. By [Iy2, Proposition 6.3(2)], ${\rm
Ext}_{\Lambda}^n(S, \Lambda)\in \mod \Lambda ^{op}$ is simple. By
[HuZ, Lemma 2.4], $S\not\subseteq
I^{0}(\Lambda)\bigoplus\cdots\bigoplus I^{n-1}(\Lambda)$. So ${\rm
Ext}_{\Lambda}^i(S, \Lambda)\cong {\rm Hom}_{\Lambda}(S,
I^i(\Lambda))=0$ for any $0\leq i \leq n-1$ by Lemma 2.2. Then from
the minimal projective resolution of $S$, we get the exact sequence:
$$0 \to P_0(S)^* \to \cdots \to P_{n-1}(S)^*\to P_n(S)^*
\to {\rm Ext}_{\Lambda}^n(S, \Lambda) \to 0$$ which is a minimal
projective resolution of ${\rm Ext}_{\Lambda}^n(S, \Lambda)$ by [M,
Proposition 4.2], where $(-)^*={\rm Hom}_{\Lambda}(-, \Lambda)$. So
$P_n(S)^* \cong P_0({\rm Ext}_{\Lambda}^n(S, \Lambda))$ is
indecomposable and hence $P_n(S)$ is also indecomposable.
$\hfill{\square}$

\vspace{0.2cm}

Denote by $\mathscr{P}^{n}(S)$ and $\mathscr{I}^{n}(S)$ the
subcategory of $\mod \Lambda$ consisting of simple modules with
projective dimension $n$ and injective dimension $n$, respectively.
Since $\mathbb{D}$ is a duality between simple $\Lambda$-modules and
simple $\Lambda ^{op}$-modules, we get easily the following result
from [Iy2, Proposition 6.3].

\vspace{0.2cm}

{\bf Lemma 3.5} {\it Let $\Lambda$ be an $(n-1)$-Auslander algebra
with $\gldim\Lambda =n$. Then the functor $\mathbb{D}\circ{\rm
Ext}_{\Lambda}^{n}(\ ,\Lambda)$ gives a bijection from
$\mathscr{P}^{n}(S)$ {\it to} $\mathscr{I}^{n}(S)$ with the inverse
${\rm Ext}_{\Lambda}^{n}(\ ,\Lambda)\circ\mathbb{D}$.}

\vspace{0.2cm}

Let $\mathscr{C}$ be a full subcategory of $\mod \Lambda$ and $n$ a
positive integer. Recall from [AuR] that $\mathscr{C}$ is said to be
{\it contravariantly finite} in $\mod\Lambda$ if for any
$M\in\mod\Lambda$, there exists a morphism $C_{M}\rightarrow M$ with
$C_{M}\in\mathscr{C}$ such that ${\rm
Hom}_{\Lambda}(C,C_{M})\rightarrow {\rm
Hom}_{\Lambda}(C,M)\rightarrow 0$ is exact for any $C\in
\mathscr{C}$. Dually, the notion of {\it covariantly finite
subcategories} of $\mod\Lambda$ is defined. A full subcategory of
$\mod\Lambda$ is said to be {\it functorially finite} in
$\mod\Lambda$ if it is both contravariantly finite and covariantly
finite in $\mod \Lambda$. We denote by $^{\bot_n}\mathscr{C}= \{
X\in \mod \Lambda \ |\ {\rm Ext}_{\Lambda}^{i}(X, C)=0$ for any $C
\in \mathscr{C}$ and $1 \leq i \leq n \}$, and
$\mathscr{C}^{\bot_n}= \{ X\in \mod \Lambda \ |\ {\rm
Ext}_{\Lambda}^{i}(C, X)=0$ for any $C \in \mathscr{C}$ and $1 \leq
i \leq n \}$.

\vspace{0.2cm}

{\bf Definition 3.6} ([Iy3]) Let $\mathscr{C}$ be a functorially
finite subcategory of $\mod\Lambda$. For a positive integer $n$,
$\mathscr{C}$ is called a {\it maximal} $n$-{\it orthogonal
subcategory} of $\mod\Lambda$ if
$\mathscr{C}={^{\bot_n}\mathscr{C}}=\mathscr{C}^{\bot_n}$.

\vspace{0.2cm}

From the definition above, we get easily that both  $\Lambda$ and
$\mathbb{D}\Lambda ^{op}$ are in any maximal $n$-orthogonal
subcategory of $\mod \Lambda$. For a module $M \in \mod \Lambda$, we
use $\add _{\Lambda}M$ to denote the subcategory of $\mod \Lambda$
consisting of all modules isomorphic to direct summands of finite
direct sums of copies of $_{\Lambda}M$. Then
$\add_{\Lambda}(\Lambda\bigoplus \mathbb{D}\Lambda ^{op})$ is
contained in any maximal $n$-orthogonal subcategory of $\mod
\Lambda$. On the other hand, it is easy to see that if
$\add_{\Lambda}(\Lambda\bigoplus \mathbb{D}\Lambda ^{op})$ is a
maximal $n$-orthogonal subcategory of $\mod \Lambda$, then
$\add_{\Lambda}(\Lambda\bigoplus \mathbb{D}\Lambda ^{op})$ is the
unique maximal $n$-orthogonal subcategory of $\mod \Lambda$. In this
case, we say that $\Lambda$ admits a {\it trivial maximal} $n$-{\it
orthogonal subcategory} of $\mod \Lambda$ (see [HuZ]).

For a positive integer $n$, we proved in [HuZ, Proposition 3.2] that
$\Lambda$ admits no maximal $j$-orthogonal subcategories of $\mod
\Lambda$ for any $j\geq n$ if $\id _{\Lambda}\Lambda=n$ (especially,
if $\gldim \Lambda =n$). Furthermore, in [HuZ] we gave an equivalent
characterization for the existence of trivial maximal
$(n-1)$-orthogonal subcategories of $\mod\Lambda$ over an
$(n-1)$-Auslander algebra $\Lambda$ with $\gldim \Lambda=n$ as
follows.

\vspace{0.2cm}

{\bf Lemma 3.7} ([HuZ, Corollary 3.10]) {\it Let $\Lambda$ be an
$(n-1)$-Auslander algebra with $\gldim \Lambda=n$. Then the
following statements are equivalent.

(1) $\Lambda$ admits a trivial maximal $(n-1)$-orthogonal
subcategory $\add_{\Lambda}(\Lambda\bigoplus \mathbb{D}\Lambda
^{op})$ of $\mod \Lambda$.

(2) A simple module $S \in \mod \Lambda$ is injective if $\pd
_{\Lambda}S=n$.}

\vspace{0.2cm}

For a positive integer $n$, recall from [FGR] that $\Lambda$ is
called {\it $n$-Gorenstein} if $\pd _{\Lambda}I^{i}(\Lambda)\leq i$
for any $0\leq i\leq n-1$. By [FGR, Theorem 3.7], the notion of
$n$-Gorenstein is left and right symmetric. Recall from [B] that
$\Lambda$ is called {\it Auslander-Gorenstein} if $\Lambda$ is
$n$-Gorenstein for all $n$ and both $\id _{\Lambda}\Lambda$ and $\id
_{\Lambda ^{op}}\Lambda$ are finite.

\vspace{0.2cm}

{\bf Lemma 3.8} {\it Assume that $\id _{\Lambda}\Lambda=\id
_{\Lambda^{op}}\Lambda =n(<\infty)$. Then we have the following

(1) ([IS, Proposition 1(1)]) $\pd _{\Lambda}X=n$ or $\infty$ for any
non-zero submodule $X$ of $I^{n}(\Lambda)$.

(2) ([IS, Corollary 7(2)]) If $\Lambda$ is Auslander-Gorenstein and
$I\in \mathscr{IP}^{n}(\Lambda)$, then $I \cong I^{0}(S)$ for some
simple module $S \in \mod \Lambda$ with $\pd _{\Lambda}S=n$ or
$\infty$.}

\vspace{0.2cm}

For a module $M \in \mod \Lambda$, the {\it grade} of $M$, denoted
by $\grade M$, is defined as ${\rm inf}\{n\geq 0\ |\ {\rm
Ext}_{\Lambda}^{n}(M,\Lambda)\neq0\}$ (see [AuB]).

\vspace{0.2cm}

{\bf Lemma 3.9} ([Iy1, Proposition 2.4]) {\it Let $\Lambda$ be
$n$-Gorenstein. Then the subcategory $\{X\in \mod \Lambda \ |\
\grade X \geq n\}$ of $\mod \Lambda$ is closed under submodules and
factor modules.}

\vspace{0.2cm}

{\bf Lemma 3.10} ([HuZ, Lemma 3.4]) {\it If $\gldim \Lambda=n\geq2$
and $\mathscr{C}$ is a subcategory of $\mod \Lambda$ such that
$\Lambda \in\mathscr{C}$ and ${\rm
Ext}_{\Lambda}^{j}(\mathscr{C},\mathscr{C})=0$ for any $1\leq j\leq
n-1$, then $\grade M=n$ for any $M\in \mathscr{C}$ without
projective direct summands.}

\vspace{0.5cm}

\centerline{\large \bf 4. The existence of trivial maximal
orthogonal subcategories}

\vspace{0.2cm}

In this section, we will mainly study the properties of
$(n-1)$-Auslander algebras with $\gldim \Lambda =n$ admitting a
trivial maximal $(n-1)$-subcategory. We will prove that for an
$(n-1)$-Auslander algebra $\Lambda$ with $\gldim \Lambda =n$, if
$\Lambda$ admits a trivial maximal $(n-1)$-orthogonal subcategory of
$\mod\Lambda$, then $\Lambda$ is a Nakayama algebra and the
projective dimension or injective dimension of any indecomposable
module in $\mod\Lambda$ is at most $n-1$. As a consequence, we have
that for an Auslander algebra $\Lambda$ with $\gldim \Lambda =2$, if
$\Lambda$ admits a trivial maximal 1-orthogonal subcategory of
$\mod\Lambda$, then $\Lambda$ is a tilted algebra of finite
representation type.

\vspace{0.2cm}

{\bf Lemma 4.1} {\it Let $\Lambda$ be an $(n-1)$-Auslander algebra
with $\gldim\Lambda=n$ and $S\in \mod \Lambda$ a simple module.

(1) If $\pd _{\Lambda}S\leq n-1$, then $I^{0}(S)$ is projective.

(2) If $\pd _{\Lambda}S=n$, then $\pd _{\Lambda}I^{0}(S)=n$.}

\vspace{0.2cm}

{\it Proof.} For any $0 \leq i \leq n$, if $\pd _{\Lambda}S=i$, then
${\rm Hom}_{\Lambda}(S, I^{i}(\Lambda))\cong {\rm
Ext}_{\Lambda}^{i}(S, \Lambda)\neq 0$. It follows that $I^{0}(S)$ is
isomorphic to a direct summand of $I^{i}(\Lambda)$. Notice that
$\Lambda$ is an $(n-1)$-Auslander algebra, then (1) follows
trivially, and (2) follows from Lemma 3.8(1). $\hfill{\square}$

\vspace{0.2cm}

For a module $M\in\mod\Lambda$, we use $\L(M)$ to denote the length
of $M$.

\vspace{0.2cm}

{\bf Lemma 4.2} {\it Let $\Lambda$ be an $(n-1)$-Auslander algebra
with $\gldim\Lambda=n$ admitting a trivial maximal
$(n-1)$-orthogonal subcategory of $\mod\Lambda$ and $M\in \mod
\Lambda$ indecomposable. If} $\L(M) \geq 2$ {\it or $M$ is not
injective, then the following equivalent conditions hold true.

(1) $\pd _{\Lambda}S\leq n-1$ for any simple submodule $S$ of $M$.

(2) $I^{0}(M)$ is projective.}

\vspace{0.2cm}

{\it Proof.} By Lemma 3.7, a simple module $S\in \mod \Lambda$ is
injective if $\pd _{\Lambda}S=n$. Because $M\in \mod \Lambda$ is
indecomposable, we have that $\pd _{\Lambda}S\leq n-1$ for any
simple submodule $S$ of $M$ and the assertion (1) holds true.
Otherwise, $M\cong S$, which contradicts with the assumption that
$\L(M) \geq 2$ or $M$ is not injective.

It suffices to prove $(1)\Rightarrow (2)$. By Lemma 4.1(1), it is
easy to get the desired conclusion. $\hfill{\square}$

\vspace{0.2cm}

The following proposition plays a crucial role in the proof of the
main result in this paper.

\vspace{0.2cm}

{\bf Proposition 4.3} {\it Let $\Lambda$ be an $(n-1)$-Auslander
algebra with $\gldim\Lambda=n$ and $0 \leq k \leq n$. If $\Lambda$
admits a trivial maximal $(n-1)$-orthogonal subcategory of
$\mod\Lambda$, then for any indecomposable non-projective-injective
module $M \in \mod \Lambda$ with $\pd_{\Lambda}M=n-k$, there exists
a simple module $S\in \mod \Lambda$ such that $\pd_{\Lambda}S=n$ and
$M\cong \Omega^{k}S$.}

 \vspace{0.2cm}

{\it Proof.} For the case $k=0$, it suffices to prove that
$\L(M)=1$. Then $M$ is simple and it is injective by Lemma 3.7. Thus
the assertion follows.

Assume that $\L(M)\geq 2$. By Lemma 4.2, $\pd _{\Lambda}S\leq n-1$
for any simple submodule $S$ of $M$ and $I^{0}(M)$ is projective.

If $M$ is injective, then $M\cong I^{0}(S)$ for some simple
$\Lambda$-module $S$ with $\pd _{\Lambda}S=n$ by Lemma 3.8(2), which
is a contradiction. Now assume that $\id _{\Lambda}M \geq 1$. By
Corollary 2.4,
$$0\rightarrow P_{n}(M)\rightarrow\cdots\rightarrow
P_{1}(M)\rightarrow P_{0}(M)\rightarrow
I^{0}(M)\stackrel{\pi}{\rightarrow} I^{0}(M)/M \rightarrow 0$$ is a
minimal projective resolution of $I^{0}(M)/M$ and $\pd
_{\Lambda}I^{0}(M)/M=n+1$, which contradicts with $\gldim\Lambda=n$.
So the case for $k=0$ is proved.

For the case $k=n$, we have that $M$ is projective. Then $M$ is not
injective by assumption. Because $\gldim\Lambda=n$, $\id
_{\Lambda}M\leq n$. On the other hand, because $\Lambda$ admits a
trivial maximal $(n-1)$-orthogonal subcategory of $\mod \Lambda$,
${\rm Ext}_{\Lambda}^{j}(\mathbb{D}\Lambda ^{op},\Lambda)=0$ for any
$1\leq j\leq n-1$. Then it is not difficult to show that
$\id_{\Lambda}M=n$. By Lemma 3.3, there exists an indecomposable
injective module $T \in \mod \Lambda$ with $\pd_{\Lambda}T=n$ such
that $M\cong \Omega^{n}T$. By the above argument, $T$ is simple.

Now assume that $1\leq k\leq n-1$. Then $\pd_{\Lambda}M=n-k\neq 0$.
We claim that $M$ is not injective. Otherwise, if $M$ is injective,
then the minimal projective resolution of $M$ splits because ${\rm
Ext}_{\Lambda}^{j}(\mathbb{D}\Lambda ^{op},\Lambda)=0$ for any
$1\leq j\leq n-1$. It follows that $M$ is projective, which is a
contradiction. The claim is proved. Then by Lemma 4.2,
$\pd_{\Lambda}S \leq n-1$ for any simple submodule $S$ of $M$ and
$I^{0}(M)$ is projective. In the following, we will prove the
assertion by induction on $k$.

If $k=1$, then $\pd_{\Lambda}M=n-1$. By Lemma 2.6 and Corollary 2.4,
$\pd_{\Lambda}I^{0}(M)/M=n$. So $I^{0}(M)/M \cong S$ for some simple
module $S$ with $\pd_{\Lambda}S=n$ by the above argument and hence
$M\cong \Omega^{1}S$.

Assume that $2\leq k\leq n-1$ and $\pd_{\Lambda}M=n-k$. By Corollary
2.4, we have a minimal projective resolution of $I^{0}(M)/M$ as
follows.
$$0\rightarrow
P_{n-k}(M)\rightarrow\cdots\rightarrow P_{1}(M)\rightarrow
P_{0}(M)\rightarrow I^{0}(M)\stackrel{\pi}{\rightarrow} I^{0}(M)/M
\rightarrow 0.$$ So $\pd_{\Lambda}I^{0}(M)/M=n-(k-1)$ and
$I^{0}(M)/M$ is indecomposable by Lemma 2.6. By the induction
hypothesis, $I^{0}(M)/M \cong \Omega^{k-1}S$ for some simple module
$S\in \mod \Lambda$ with $\pd_{\Lambda}S=n$. It follows that $M\cong
\Omega^{k}S$. $\hfill{\square}$

\vspace{0.2cm}

As an application of Proposition 4.3, we get the following theorem,
which is one of the main results in this section.

\vspace{0.2cm}

{\bf Theorem 4.4} {\it Let $\Lambda$ be an $(n-1)$-Auslander algebra
with $\gldim\Lambda=n$. If $\Lambda$ admits a trivial maximal
$(n-1)$-orthogonal subcategory of $\mod\Lambda$, then

(1) $\pd _{\Lambda}M+\id _{\Lambda}M=n$ for any
non-projective-injective indecomposable module $M \in \mod
\Lambda$.}

(2) $\pd_{\Lambda}M\leq n-1$ or $\id_{\Lambda}M\leq n-1$ for any
indecomposable module $M\in\mod\Lambda$.

\vspace{0.2cm}

{\it Proof.} (1) Assume that $\pd _{\Lambda}M=n-k$ for some $0\leq k
\leq n$. By Proposition 4.3, there exists a simple module $S\in \mod
\Lambda$ such that $\pd_{\Lambda}S=n$ and $M\cong \Omega^{k}S$. Then
$S$ is injective by Lemma 3.4 and so $S$ is isomorphic to a direct
summand of $\mathbb{D}\Lambda$. Because $\Lambda$ is an
$(n-1)$-Auslander algebra, $P_{i}(S)$ is injective for any $0 \leq i
\leq n-1$ by Lemma 3.2. Then, by Lemma 2.6, the following exact
sequence:
$$0\rightarrow \Omega ^{k}S \rightarrow P_{k-1}(S)\rightarrow\cdots\rightarrow
P_{1}(S)\rightarrow P_{0}(S)\rightarrow S \rightarrow 0$$ is a
minimal injective resolution of $\Omega ^{k}S(\cong M)$ and $\id
_{\Lambda}M=k$.

(2) follows from (1) immediately.  $\hfill{\square}$

\vspace{0.2cm}

The following result is another application of Proposition 4.3.

\vspace{0.2cm}

{\bf Proposition 4.5} {\it Let $\Lambda$ be an $(n-1)$-Auslander
algebra with $\gldim\Lambda=n$ and $S\in \mod\Lambda$ a simple
module with $\pd_{\Lambda}S=n$. If $\Lambda$ admits a trivial
maximal $(n-1)$-orthogonal subcategory of $\mod\Lambda$, then
$\Omega^{i}S$ is simple and $P_{i}(S)$ is indecomposable for any
$0\leq i\leq n$.}

\vspace{0.2cm}

{\it Proof.} Assume that $S\in \mod\Lambda$ is a simple module with
$\pd_{\Lambda}S=n$. By Lemma 3.7, $S$ is injective. It follows from
Lemma 3.2(2) that the minimal projective resolution of $S$:
$$0\rightarrow P_{n}(S)\rightarrow\cdots\rightarrow
P_{1}(S)\rightarrow P_{0}(S)\rightarrow S\rightarrow 0$$ is a
minimal injective resolution of $P_{n}(S)$.

We proceed by induction on $i$. The case for $i=0$ holds true
trivially, and the case for $i=n$ follows from Lemma 3.4 and the
dual version of Proposition 4.3.

Now assume that $1\leq i \leq n-1$ and $S^{'}\in \mod \Lambda$ is a
simple submodule of $\Omega^{i}S$. Because $\Lambda$ is an
$(n-1)$-Auslander algebra and $S$ is injective, $P_{0}(S)$ is
projective-injective and indecomposable. So $S^{'}$ is the unique
simple submodule of $P_{0}(S)$ and hence $I^{0}(S^{'})=P^{0}(S)$. By
Lemma 2.2, ${\rm Ext}_{\Lambda}^{n-1}(S^{'},P_{n}(S))\cong {\rm
Hom}_{\Lambda}(S^{'},\Omega^{1}S)\not=0$, which implies that
$\pd_{\Lambda}S^{'}\geq n-1$. Because $\gldim\Lambda=n$, it is easy
to see that $\pd_{\Lambda}S^{'}= n-1$. Then by Theorem 4.4,
$\id_{\Lambda}S^{'}=1$.

Connecting a minimal projective resolution  and a minimal injective
resolution of $S^{'}$, then by Lemma 2.6, the following exact
sequence is a minimal projective resolution of $I^{1}(S^{'})$:
$$0\rightarrow P_{n-1}(S^{'})\rightarrow\cdots\rightarrow
P_{0}(S^{'})\rightarrow I^{0}(S^{'})(\cong P_{0}(S))\rightarrow
I^{1}(S^{'})\rightarrow0$$ with $I^{1}(S^{'})$ indecomposable. So
$\pd _{\Lambda}I^{1}(S^{'})=n$ and hence $I^{1}(S^{'})$ is simple by
Lemma 3.2(1). It follows that $S\cong I^{1}(S^{'})$ and then
$\Omega^{1}S \cong S^{'}$ is simple. Thus $P_{1}(S)$ is
indecomposable. The case for $i=1$ is proved.

Assume that $2\leq i\leq n-1$. By Lemma 2.2, ${\rm
Ext}_{\Lambda}^{n-i}(S^{'},P_{n}(S))\cong {\rm
Hom}_{\Lambda}(S^{'},\Omega^{i}S)\not=0$. So
$\pd_{\Lambda}S^{'}(=t)\geq n-i$.

Consider the following commutative diagram with exact rows:
$$\xymatrix{& 0 \ar[r] & S^{'}\ar[r]
\ar[d]^{\alpha} & P_{i-1}(S) \ar[r]^{\pi}\ar@{=}[d] & M \ar[r] \ar[d]^{\beta}& 0\\
& 0 \ar[r] & \Omega^{i}S\ar[r] & P_{i-1}(S) \ar[r] & \Omega^{i-1}S
\ar[r] & 0}$$ where $M=P_{i-1}(S)/S^{'}$, $\alpha$ is an embedding
homomorphism and $\beta$ is an induced homomorphism. By the
induction hypothesis, $\Omega^{i-1}S$ is simple and hence
$P_{i-1}(S)$ is indecomposable. Then, by Lemma 2.6, $M$ is
indecomposable and $\pi$ is right minimal. It follows that $\pd
_{\Lambda}M=t+1$. Thus $M\cong \Omega^{n-t-1}S^{''}$ for some simple
module $S^{''}\in \mod \Lambda$ with $\pd_{\Lambda}S^{''}=n$ by
Proposition 4.3. Because $i\geq n-t-1$, $M$ is simple by the
induction hypothesis. It is clear that $\beta$ is an epimorphism and
so it is an isomorphism, which implies that $\alpha$ is an
isomorphism and $\Omega^{i}S\cong S^{'}$ is simple. It follows that
$P_i(S)$ is indecomposable. $\hfill{\square}$

\vspace{0.2cm}

As an application of Proposition 4.5, we get the following

\vspace{0.2cm}

{\bf Corollary 4.6} {\it Let $\Lambda$ be an $(n-1)$-Auslander
algebra with $\gldim\Lambda=n$. If $\Lambda$ admits a trivial
maximal $(n-1)$-orthogonal subcategory of $\mod\Lambda$, then for
any $0\leq i\leq n$, $\Omega^{i}$ gives a bijection from $\{[S]\ |\
S\in \mod \Lambda$ is simple with $\pd_{\Lambda}S=n\}$ to $\{[S]\ |\
S\in \mod \Lambda$ is non-projective-injective and simple with
$\pd_{\Lambda}S=n-i\}$, where $[S]$ is the isomorphic class of
modules in $\mod \Lambda$ containing $S$.}

\vspace{0.2cm}

{\it Proof.} By Proposition 4.5, $\Omega^{i}: \{[S]\ |\ S\in \mod
\Lambda$ is simple with $\pd_{\Lambda}S=n\}\to \{[S]\ |\ S\in \mod
\Lambda$ is simple with $\pd_{\Lambda}S=n-i\}$ is a map. By
Proposition 4.3, $\Omega^{i}$ is epic. On the other hand, any simple
module $S\in \mod \Lambda$ with $\pd _{\Lambda}S=n$ is injective by
Lemmas 3.7 or Theorem 4.4, so the minimal projective resolution of
$S$ is a minimal injective resolution of $P_n(S)$ by Lemma 3.2(2).
In particular, $0 \to \Omega ^{i}S \to P_{i-1}(S) \to \cdots \to
P_{1}(S) \to P_{0}(S) \to S \to 0$ is a minimal injective resolution
of $\Omega ^{i}S$. Then it is easy to see that $\Omega ^i$ is monic.
$\hfill{\square}$

\vspace{0.2cm}

We give another application of Proposition 4.5 as follows.

\vspace{0.2cm}

{\bf Corollary 4.7}  {\it Let $\Lambda$ be an $(n-1)$-Auslander
algebra with $\gldim\Lambda=n$. If $\Lambda$ admits a trivial
maximal $(n-1)$-orthogonal subcategory of $\mod\Lambda$, then for
any indecomposable projective module $P\in \mod \Lambda$, either $P$
or the radical $\rad P$ of $P$ is simple.}

\vspace{0.2cm}

{\it Proof.} Let $P \in \mod \Lambda$ be an indecomposable
projective module. Then there exists a unique (up to isomorphisms)
simple module $S \in \mod \Lambda$ such that $P\cong P_0(S)$. If $S$
is projective, then $P$ is simple. Now suppose $\pd
_{\Lambda}S=n-k>0$. Then by Proposition 4.3, there exists a simple
module $S^{'}\in \mod \Lambda$ with $\pd _{\Lambda}S^{'}=n$ such
that $S \cong \Omega ^{k}S^{'}$. By Proposition 4.5, $\rad P(\cong
\Omega ^{1}S \cong \Omega ^{k+1}S^{'})$ is simple. $\hfill{\square}$

\vspace{0.2cm}

{\bf Definition 4.8} ([AuRS]) $\Lambda$ is called a {\it Nakayama
algebra} if every indecomposable projective module and every
indecomposable injective module in $\mod \Lambda$ have a unique
composition series.

\vspace{0.2cm}

The following theorem is another main result in this section.

\vspace{0.2cm}

{\bf Theorem 4.9}  {\it Let $\Lambda$ be an $(n-1)$-Auslander
algebra with $\gldim\Lambda=n$. If $\Lambda$ admits a trivial
maximal $(n-1)$-orthogonal subcategory of $\mod\Lambda$, then
$\Lambda$ is a Nakayama algebra.}

\vspace{0.2cm}

{\it Proof.} Let $P \in \mod \Lambda$ be an indecomposable
projective module. Then $\rad P$ is the unique maximal submodule of
$P$. It follows from Corollary 4.7 that $P$ has a unique composition
series of length at most two.

Note that $\Lambda^{op}$ is also an $(n-1)$-Auslander algebra
admitting a trivial maximal $(n-1)$-orthogonal subcategory of
$\mod\Lambda^{op}$. So by the above argument, every indecomposable
projective module in $\mod \Lambda ^{op}$ has a unique composition
series of length at most two. Applying the functor $\mathbb{D}$,
then we get that every indecomposable injective module in
$\mod\Lambda$ has a unique composition series of length at most two.
Thus $\Lambda$ is a Nakayama algebra. $\hfill{\square}$

\vspace{0.2cm}

The following example illustrates that there exists a basic and
connected $(n-1)$-Auslander algebra $\Lambda$ with $\gldim \Lambda
=n$, which is a Nakayama algebra, but admits no maximal
$(n-1)$-orthogonal subcategories of $\mod \Lambda$. It means that
the converse of Theorem 4.9 does not hold true in general.

\vspace{0.2cm}

{\bf Example 4.10} Let $\Lambda$ be a finite-dimensional algebra
over an algebraically closed field given by the quiver:

{\tiny$$\xymatrix{1 & \ar[l]_{\beta _{1}} 2 & \ar[l]_{\beta _{2}} 3
& \ar[l]_{\beta _{3}} \cdots & \ar[l]_{\beta _{n-1}} n &
\ar[l]_{\beta _{n}} n+1 & \ar[l]_{\beta _{n+1}} n+2 & \ar[l]_{\beta
_{n+2}} \cdots & \ar[l]_{\beta _{2n-2}} 2n-1 & \ar[l]_{\beta
_{2n-1}} 2n & \ar[l]_{\beta _{2n}} 2n+1}$$}modulo the ideal
generated by $\{\beta_{i}\beta_{i+1}\ |\ 1 \leq i \leq 2n-1$ but
$i\neq n\}$. Then $\Lambda$ is a basic and connected
$(n-1)$-Auslander algebra with $\gldim \Lambda =n$. By [AsSS,
Chapter V, Theorem 3.2] (see Lemma 5.2 below), $\Lambda$ is a
Nakayama algebra. We use $P(i)$, $I(i)$ and $S(i)$ to denote the
projective, injective and simple modules corresponding to the vertex
$i$ for any $1 \leq i \leq 2n+1$. Because $P(n+2)=I(n)$ is not
simple, it follows from [AsSS, Chapter IV, Proposition 3.11] that
$0\to P(n+1) \to S(n+1)\bigoplus P(n+2) \to I(n+1) \to 0$ is an
almost split sequence. So ${\rm Ext}_{\Lambda}^1(I(n+1), P(n+1))\neq
0$ and hence there does not exist a maximal $j$-orthogonal
subcategory of $\mod\Lambda$ for any $j \geq 1$.

\vspace{0.2cm}

In the rest of this section, we will apply Theorems 4.4 and 4.9 to
Auslander algebras with global dimension 2. As a result, we can give
a connection between Auslander algebras and tilted algebras. We
first recall some notions from [HRS] and [HRi].

\vspace{0.2cm}

{\bf Definition 4.11} (1) ([HRS]) $\Lambda$ is called {\it almost
hereditary} if the following conditions are satisfied: (a) $\gldim
\Lambda\leq 2$; and (b) If $X \in \mod \Lambda$ is indecomposable,
then either $\pd _{\Lambda}X\leq 1$ or $\id _{\Lambda}X\leq 1$.

(2) ([HRS]) $\Lambda$ is called {\it quasi-tilted} if $\Lambda ={\rm
End}(T)^{op}$, where $\mathcal{H}$ is a locally finite hereditary
abelian $R$-category and $T$ is a tilting object in $\mathcal{H}$.

(3) ([HRi]) $\Lambda$ is called {\it tilted} if $\Lambda$ is of the
form $\Lambda ={\rm End}(T_{\Gamma})$, where $T_{\Gamma}$ is a
tilting module and $\Gamma$ is a hereditary Artinian algebra. It is
trivial that a tilted algebra is quasi-tilted.

\vspace{0.2cm}

Now we are in a position to give the following result.

\vspace{0.2cm}

{\bf Corollary 4.12} {\it Let $\Lambda$ be an Auslander algebra with
$\gldim\Lambda=2$. If $\Lambda$ admits a trivial maximal
1-orthogonal subcategory of $\mod\Lambda$, then $\Lambda$ is a
tilted algebra of finite representation type.}

\vspace{0.2cm}

{\it Proof.} Let $\Lambda$ be an Auslander algebra with
$\gldim\Lambda=2$ admitting a trivial maximal 1-orthogonal
subcategory of $\mod\Lambda$. Then $\Lambda$ is an almost hereditary
algebra of finite representation type by Theorems 4.4 and 4.9. So
$\Lambda$ is quasi-tilted by [HRS, Chapter III, Theorem 2.3] and
hence it is tilted by [HRS, Chapter III, Corollary 3.6].
$\hfill{\square}$

\vspace{0.2cm}

{\it Remark.} (1) Let $\Lambda$ be an Auslander algebra (of finite
representation type) with $\gldim\Lambda=2$. If $\Lambda$ admits a
non-trivial maximal 1-orthogonal subcategory of $\mod\Lambda$ (note:
Iyama in [Iy5] constructed an example to illustrate that this may
occur), then $\Lambda$ is not almost hereditary because any maximal
1-orthogonal subcategory (if it exists) for an almost hereditary
algebra is trivial by [HuZ, Theorem 3.15]. So $\Lambda$ is not
(quasi-)tilted by [HRS, Chapter III, Theorem 2.3].

(2) In the statement of Corollary 4.12, the conditions ``$\Lambda$
is an Auslander algebra" and ``$\Lambda$ is a tilted algebra of
finite representation type" cannot be exchanged. For example, let
$\Lambda$ be a finite-dimensional algebra given by the quiver:
$$\xymatrix{1&2\ar[l]_{\alpha_{1}}&3\ar[l]_{\alpha_{2}}
&4\ar[l]_{\alpha_{3}}&5\ar[l]_{\alpha_{4}}}$$ modulo the ideal
generated by $\{\alpha_{1}\alpha_{2}\alpha_{3}\alpha_{4}\}$. Then
$\Lambda$ is a tilted algebra of finite representation type (cf.
[AsSS, p.323]), and $\Lambda$ admits a trivial maximal 1-orthogonal
subcategory $\add _{\Lambda} \bigoplus_{i=1}^5P(i)\bigoplus
I(3)\bigoplus I(4)\bigoplus I(5)$ of $\mod\Lambda$. However,
$\Lambda$ is not an Auslander algebra because $\pd
_{\Lambda}I^1(\Lambda)=2$.

\vspace{0.5cm}

\centerline{\Large \bf 5. The case for finite-dimensional algebras}

\vspace{0.2cm}

In this section, $\Lambda$ is a finite-dimensional algebra over an
algebraically closed field $K$. As an application of the results
obtained in Section 4, we will give an explicit description of the
quiver of an $(n-1)$-Auslander algebra $\Lambda$ with $\gldim
\Lambda =n$ admitting a trivial maximal $(n-1)$-orthogonal
subcategory of $\mod \Lambda$. We begin with some preliminary
results.

\vspace{0.2cm}

{\bf Lemma 5.1} {\it Let $\{P_{1},P_{2},\cdots ,P_{m}\}$ be a
complete set of non-isomorphic indecomposable projective modules in
$\mod \Lambda$. Then $\Lambda$ is connected if and only if there
does not exist a non-trivial partition $J_{1}\dot{\bigcup}J_{2}$ of
the set $\{P_{1},P_{2},\cdots ,P_{m}\}$ such that $P_i\in J_{1}$ and
$P_j\in J_{2}$ imply ${\rm Hom}_{\Lambda}(P_i, P_j)=0={\rm
Hom}_{\Lambda}(P_j, P_i)$.}

\vspace{0.2cm}

{\it Proof.} By [AsSS, Chapter II, Lemma 1.6].  $\hfill{\square}$

\vspace{0.2cm}

{\bf Lemma 5.2} ([AsSS, Chapter V, Theorem 3.2]) {\it Let $\Lambda$
be basic and connected. Then $\Lambda$ is a Nakayama algebra if and
only if its ordinary quiver is one of the following two quivers:

(1) $1\leftarrow 2\leftarrow3\leftarrow\cdots \leftarrow
n-1\leftarrow n$;

(2) $$
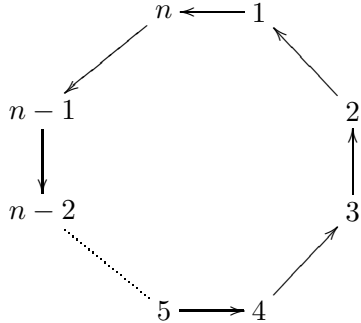
\xymatrix{
&n\ar[ld]&1\ar[l]&\\
n-1\ar[d]&&&2\ar[ul]\\
n-2&&&3\ar[u]\\
&5\ar@{.}[lu]\ar[r]&4\ar[ur]& }
$$ (with $n\geq 1$ points).}

\vspace{0.2cm}

The following proposition is useful for proving the main result in
this section.

\vspace{0.2cm}

{\bf Proposition 5.3} {\it Let $\Lambda$ be a basic and connected
$(n-1)$-Auslander algebra with $\gldim\Lambda=n$. If $\Lambda$
admits a trivial maximal $(n-1)$-orthogonal subcategory of
$\mod\Lambda$, then there exists a unique simple module
$S\in\mod\Lambda$ with $\pd_{\Lambda}S=n$.}

\vspace{0.2cm}

{\it Proof.} Since $\Lambda$ is connected, it is not difficult to
see that there does not exist a simple projective-injective module
in $\mod \Lambda$ by Lemma 5.1. Then by Corollary 4.6, every simple
module in $\mod \Lambda$ is of the form $\Omega^{i}S$ for some
simple module $S\in\mod\Lambda$ with $\pd_{\Lambda}S=n$ and $0\leq
i\leq n$. Because $\gldim \Lambda =n$, there exists a simple module
$S\in \mod \Lambda$ with $\pd _{\Lambda}S=n$.

Assume that $\{S_{1},S_{2},\cdots,S_{t} \}$ is a complete set of
non-isomorphic simple modules in $\mod \Lambda$ with projective
dimension $n$. It suffices to prove $t=1$. Suppose $t \geq 2$. By
Corollary 4.6, we get a complete set of non-isomorphic
indecomposable projective modules in $\mod \Lambda$ as follows.
$$\big\{P_{0}(S_{1}), P_{1}(S_{1})(=P_{0}(\Omega^{1}{S_{1}})), P_{2}(S_{1})(=P_{0}(\Omega^{2}{S_{1}})),
\cdots, P_{n}(S_{1})(=P_{0}(\Omega^{n}{S_{1}}))$$
$$P_{0}(S_{2}), P_{1}(S_{2})(=P_{0}(\Omega^{1}{S_{2}})), P_{2}(S_{2})(=P_{0}(\Omega^{2}{S_{2}})),
\cdots, P_{n}(S_{2})(=P_{0}(\Omega^{n}{S_{2}}))$$
$$\cdots\cdots\cdots\cdots$$
$$P_{0}(S_{t}), P_{1}(S_{t})(=P_{0}(\Omega^{1}{S_{t}})), P_{2}(S_{t})(=P_{0}(\Omega^{2}{S_{t}})),
\cdots, P_{n}(S_{t})(=P_{0}(\Omega^{n}{S_{t}}))\big\}.$$

In the following, we will show $${\rm
Hom}_{\Lambda}(P_{i}(S_{j}),P_{k}(S_{l}))=0={\rm
Hom}_{\Lambda}(P_{k}(S_{l}),P_{i}(S_{j})) \eqno{(*)}$$ for any
$0\leq i,\ k\leq n$ and $1\leq j\neq l\leq t$.

Notice that $S_i$ is injective for any $1 \leq i \leq t$ by Lemma
3.7, then it follows from Lemma 3.2(2) that the minimal projective
resolution of $S_l$:
$$0\to P_n(S_l)\to \cdots \to P_1(S_l)\to P_0(S_l) \to S_l \to
0 \eqno{(**)}$$ is a minimal injective resolution of $P_n(S_l)$. We
split $(**)$ to the following $n$ short exact sequences:
$$0\rightarrow P_{n}(S_{l})\rightarrow P_{n-1}(S_{l})\rightarrow
\Omega^{n-1}S_{l}\rightarrow 0  \eqno{(1)}$$
$$\cdots\cdots\cdots\cdots$$
$$0\rightarrow \Omega^{2}S_{l}\rightarrow P_{1}(S_{l})\rightarrow
\Omega^{1}S_{l}\rightarrow 0  \eqno{(n-1)}$$
$$0\rightarrow \Omega^{1}S_{l}\rightarrow P_{0}(S_{l})\rightarrow
S_{l}\rightarrow 0  \eqno{(n)}$$ Since $S_{j}\not\cong S_{l}$ and
$\Omega^{i}S_{j}$ and $\Omega^{i}S_{l}$ are simple for any $0\leq i
\leq n$ by Proposition 4.5, we get that $\Omega^{i}S_{j}\not
\cong\Omega^{i}S_{l}$ for any $0\leq i \leq n$ by Corollary 4.6
(note: $\Omega^{n}S_{j}\cong P_n(S_j)$ and $\Omega^{n}S_{l}\cong
P_n(S_l)$). So $\Omega^{i}S_{j}\not \cong\Omega^{k}S_{l}$ for any
$0\leq i, k \leq n$.

By applying the functor ${\rm Hom}_{\Lambda}(P_{i}(S_{j}),-)$ (where
$0\leq i \leq n$) to the exact sequences $(1), \cdots, (n-1), (n)$,
then we get the following exact sequences:
$$0={\rm
Hom}_{\Lambda}(P_{i}(S_{j}),\Omega ^{k+1}S_{l}) \to {\rm
Hom}_{\Lambda}(P_{i}(S_{j}),P_{k}(S_{l})) \to {\rm
Hom}_{\Lambda}(P_{i}(S_{j}),\Omega ^{k}S_{l})=0$$ for $k=0, 1,
\cdots, n-1$. So ${\rm Hom}_{\Lambda}(P_{i}(S_{j}),P_{k}(S_{l}))=0$
for any $0\leq i\leq n-1$. In addition, it is trivial that ${\rm
Hom}_{\Lambda}(P_{i}(S_{j}),P_{n}(S_{l}))=0$. This proves the left
equality of $(*)$. Dually, we get the right equality of $(*)$. Thus
we get a non-trivial partition $J_{1}\dot{\bigcup}J_{2}$ of the set
of non-isomorphic indecomposable projective modules in $\mod
\Lambda$, where $J_{1}=\big\{P_{0}(S_{1}), P_{1}(S_{1}), \cdots,
P_{n}(S_{1})\big\}$ and $J_{2}=\big\{P_{0}(S_{2}), P_{1}(S_{2}),
\cdots, P_{n}(S_{2}), \cdots, P_{0}(S_{t}), P_{1}(S_{t}), \cdots,
P_{n}(S_{t})\big\}$, which is a contradiction by Lemma 5.1. The
proof is finished. $\hfill{\square}$

\vspace{0.2cm}

Now we are in a position to give the main result in this section.

\vspace{0.2cm}

{\bf Theorem 5.4} {\it $\Lambda$ is a basic and connected
$(n-1)$-Auslander algebra with $\gldim \Lambda =n$ admitting a
trivial maximal $(n-1)$-orthogonal subcategory of $\mod\Lambda$ if
and only if $\Lambda$ is given by the quiver:
$$\xymatrix{1 & \ar[l]_{\beta _{1}} 2 & \ar[l]_{\beta _{2}} 3
& \ar[l]_{\beta _{3}} \cdots & \ar[l]_{\beta _{n}} n+1}
$$ modulo the ideal generated by $\{
\beta_{i}\beta_{i+1}| 1\leq i\leq n-1 \}$.}

\vspace{0.2cm}

{\it Proof.} We first prove the sufficiency. It is straightforward
to verify that $\Lambda$ is an $(n-1)$-Auslander algebra with
$\gldim \Lambda=n$ and admits a maximal $(n-1)$-orthogonal
subcategory $\mathscr{C}=\add_{\Lambda}\big(P(1)\bigoplus
P(2)\bigoplus P(3)\bigoplus\dots\bigoplus P(n+1)\bigoplus
S(n+1)\big)$.

We then prove the necessity. By Theorem 4.9, $\Lambda$ is a Nakayama
algebra. By Proposition 5.3 and Corollary 4.6, there exist exactly
$n+1$ non-isomorphic simple modules in $\mod\Lambda$. Because
$\gldim \Lambda =n$, there exists a simple module $S\in \mod
\Lambda$ with $\pd _{\Lambda}S=n$. By Lemma 3.7, $S$ is injective.
That is, there exists a simple injective module in $\mod\Lambda$.
Then by Lemma 5.2, the ordinary quiver of $\Lambda$ is given by
$$\xymatrix{1 & \ar[l]_{\beta _{1}} 2 & \ar[l]_{\beta _{2}} 3
& \ar[l]_{\beta _{3}} \cdots & \ar[l]_{\beta _{n}} n+1}.$$

We claim that $\beta_{i}\beta_{i+1}=0$ for any $1\leq i\leq n-1$.
Otherwise, if $\beta_{i}\beta_{i+1}\neq 0$ for some $1\leq i\leq
n-1$, then neither $P(i+2)$ nor $\rad P(i+2)(\cong \Omega ^1S(i+2))$
are simple, which is a contradiction by Corollary 4.7. The claim is
proved. Because $\Lambda$ is connected, the ideal generated by $\{
\beta_{i}\beta_{i+1}| 1\leq i\leq n-1 \}$ is exactly the non-zero
admissible ideal of $KQ$, and the assertion follows.
$\hfill{\square}$

\vspace{0.2cm}

In the following, as an application of Theorem 5.4, we will
establish the structure theorem of an $(n-1)$-Auslander algebra
$\Lambda$ admitting a trivial maximal $(n-1)$-orthogonal subcategory
of $\mod\Lambda$.

\vspace{0.2cm}

{\bf Definition 5.5} ([AuRS]) Let
$\Lambda=\Lambda_{1}\times\cdots\times\Lambda_{n}$ be a product of
indecomposable algebras and $1=e_{1}+\cdots+e_{n}$ the corresponding
decomposition of the identity element $1$ of $\Lambda$. Then the
$\Lambda _i$ are called the {\it blocks} of $\Lambda$.

\vspace{0.2cm}

{\bf Lemma 5.6} ([A, Chapter IV, Proposition 3]) {\it Let $S$ and
$T$ be simple modules in $\mod \Lambda$. Then the following
statements are equivalent.

(1) $S$ and $T$ lie in the same block.

(2) There exist simple modules $S=S_{1},S_{2},\cdots, S_{m}=T$ in
$\mod \Lambda$ such that for any $1\leq i\leq m-1$, $S_{i}$ and
$S_{i+1}$ are composition factors of an indecomposable projective
module in $\mod \Lambda$.

(3) There exist simple modules $S=S_{1},S_{2},\cdots, S_{m}=T$ in
$\mod \Lambda$ such that for any $1\leq i\leq m-1$, either
$S_{i}=S_{i+1}$ or there exists a non-split extension of one of them
by the other.}

\vspace{0.2cm}

We give some elementary properties of modules in a block as follows.

\vspace{0.2cm}

{\bf Lemma 5.7} {\it Let $\Lambda _i$ be a block of $\Lambda$ and
$M\in \mod \Lambda _i$. Then we have

(1) The submodules, quotient modules and finite direct sums of $M$
are also in $\mod \Lambda _i$.

(2) $M$ is simple in $\mod \Lambda _i$ if and only if it is simple
in $\mod \Lambda$.

(3) $M$ is projective (resp. injective) in $\mod \Lambda _i$ if and
only if it is projective (resp. injective) in $\mod \Lambda$.

(4) $P_0(M)$ (resp. $I^0(M)$)$\ \in \mod \Lambda _i$.

(5) If $\Lambda _j$ is a block of $\Lambda$ with $j\neq i$, then
${\rm Ext}_{\Lambda}^t(N, M)=0$ for any $N\in \mod \Lambda _j$ and
$t \geq 0$.}

\vspace{0.2cm}

{\it Proof.} (1) It follows from [A, p.93].

(2) If $M$ is simple in $\mod \Lambda _i$, then by (1), the
submodules of $M$ as $\Lambda$-modules are in $\mod \Lambda _i$,
which implies that $M$ is simple in $\mod \Lambda$. The converse is
trivial.

(3) If $M$ is projective in $\mod \Lambda$, then $M\bigoplus
P\cong\Lambda ^{(s)}$ for some projective module $P \in \mod
\Lambda$ and a positive integer $s$. It follows that $M\bigoplus
\Lambda _iP\cong \Lambda _iM\bigoplus \Lambda _iP\cong\Lambda
_i\Lambda ^{(s)}\cong \Lambda _i^{(s)}$ and $M$ is projective in
$\mod \Lambda _i$. The converse is trivial. By the usual duality
$\mathbb{D}$, we get another assertion.

(4) From the exact sequence $P_0(M) \to M \to 0$ in $\mod \Lambda$,
we get an exact sequence $\Lambda _iP_0(M) \to \Lambda _iM(\cong M)
\to 0$ in $\mod \Lambda _i$. Notice that $P_0(M)$ is a projective
cover of $M$ and $0\neq \Lambda _iP_0(M)$ is isomorphic to a direct
summand of $P_0(M)$, so we have that $P_0(M)\cong\Lambda _iP_0(M)
\in \mod \Lambda _i$. By the usual duality $\mathbb{D}$, we get that
$I^0(M)\in \mod \Lambda _i$.

(5) If $\Lambda _j$ is a block of $\Lambda$ with $j\neq i$, then it
follows from [A, p.93] that ${\rm Hom}_{\Lambda}(N, M)=0$ for any
$N\in \mod \Lambda _j$ and $M\in \mod \Lambda _i$. From this fact
the case for $t=0$ follows. Then by applying the functor ${\rm
Hom}_{\Lambda}(-, M)$ to a minimal projective resolution of $N$ in
$\mod \Lambda _j$, we get the assertion inductively.
$\hfill{\square}$

\vspace{0.2cm}

By the following lemma, we can draw the Auslander-Reiten quiver of
an $(n-1)$-Auslander algebra $\Lambda$ admitting a trivial maximal
$(n-1)$-orthogonal subcategory of $\mod\Lambda$.

\vspace{0.2cm}

{\bf Lemma 5.8} {\it Let $\Lambda$ be an $(n-1)$-Auslander algebra
with $\gldim\Lambda=n$ and $S\in\mod\Lambda$ a simple module with
$\pd_{\Lambda}S=n$. If $\Lambda$ admits a trivial maximal
$(n-1)$-orthogonal subcategory of $\mod\Lambda$, then
$$0\rightarrow \Omega^{i+1}S\rightarrow
P_{i-1}(S)\stackrel{\pi_{i}}{\rightarrow} \Omega^{i}S\rightarrow 0
\eqno{(*i)}$$ is an almost split sequence for any $0\leq i\leq
n-1$.}

\vspace{0.1cm}

{\it Proof.} It is obvious that $(*i)$ does not split. It suffices
to prove that there exists a homomorphism $g:\ M\rightarrow
P_{i-1}(S)$ such that $\pi_{i}g=f$ whenever $f:\ M\rightarrow
\Omega^{i}S$ is not a split epimorphism with $M\in\mod\Lambda$
indecomposable.

If $M$ is projective, then the assertion follows trivially. If
$1\leq\pd_{\Lambda}M=j\leq n$, then by Proposition 4.3, $M\cong
\Omega^{n-j}S^{'}$ for some simple module $S^{'}\in\mod\Lambda$ with
$\pd_{\Lambda}S^{'}=n$. So $M(\cong \Omega^{n-j}S^{'})$ and
$\Omega^{i}S$ are simple by Proposition 4.5. Because $f:
M\rightarrow \Omega^{i}S$ is not a split epimorphism, $f=0$. Thus
$g=0$ is desired. $\hfill{\square}$

\vspace{0.1cm}

Now we state the structure theorem of an $(n-1)$-Auslander algebra
$\Lambda$ admitting a trivial maximal $(n-1)$-orthogonal subcategory
of $\mod\Lambda$ as follows.

\vspace{0.1cm}

{\bf Theorem 5.9} {\it $\Lambda$ is an $(n-1)$-Auslander algebra
$\Lambda$ with $\gldim\Lambda=n$ admitting a trivial maximal
$(n-1)$-orthogonal subcategory of $\mod\Lambda$ if and only if it is
a finite direct product of $K$ and the algebra given by the quiver
$$\xymatrix{1 & \ar[l]_{\beta _{1}} 2 & \ar[l]_{\beta _{2}} 3
& \ar[l]_{\beta _{3}} \cdots & \ar[l]_{\beta _{n}} n+1}
$$ modulo the ideal generated by $\{
\beta_{i}\beta_{i+1}| 1\leq i\leq n-1 \}$}.

\vspace{0.1cm}

{\it Proof.} The sufficiency follows from Theorem 5.4 and Lemma 5.7.
In the following, we will prove the necessity.

Let $\{S_{1},S_{2},\cdots,S_{t} \}(t\geq1)$ be a complete set of
non-isomorphic simple modules in $\mod \Lambda$ with projective
dimension $n$. By Corollary 4.7, it is easy to see that for any
$S_i, S_j \in \{S_{1},S_{2},\cdots,S_{t} \}$ with $i\neq j$, there
does not exists an indecomposable projective module $P\in\mod
\Lambda$ such that $S_i$ and $S_j$ are composition factors of $P$.
It follows from Lemma 5.6 that $S_i$ and $S_j$ do not lie in the
same block of $\Lambda$. Then we may assume that $\Lambda = K \times
\cdots \times K \times \Lambda_{1}\times\cdots\times\Lambda_{t}$ is
a decomposition of blocks of $\Lambda$ with $S_i \in \Lambda _i$ for
any $S_i\in \{S_{1},S_{2},\cdots,S_{t} \}$. By Proposition 4.3 and
Lemma 5.8, for any $S_{i}\in$ $\{S_{1},S_{2},\cdots,S_{t} \}$, we
get a connected component of the Auslander-Reiten quiver of
$\Lambda$ as follows.
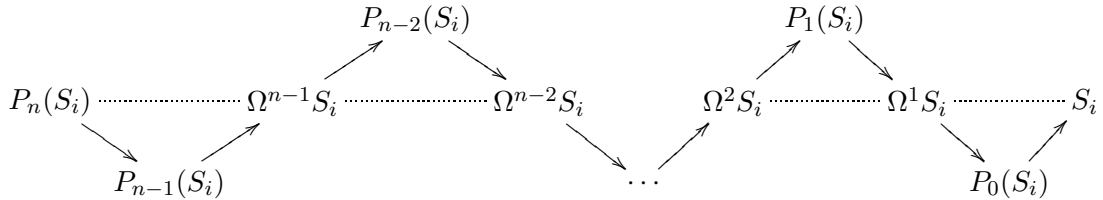
\begin{center}$\xymatrix@C=0.08cm@R=0.5cm{&&&P_{n-2}(S_i)\ar[dr]&&&&P_1(S_i)\ar[dr]\\
P_n(S_i)\ar[dr]\ar@{.}[rr]&&\Omega^{n-1}S_i\ar[ur]\ar@{.}[rr]&&\Omega^{n-2}S_i
\ar[dr]&&\Omega^2S_i\ar[ur]\ar@{.}[rr]&&\Omega^1S_i\ar[dr]\ar@{.}[rr]&&S_i\\
&P_{n-1}(S_i)\ar[ur]&&&&\cdots\ar[ur]&&&&P_0(S_i)\ar[ur]\\&&&&&(**i)}$\end{center}
By Lemmas 5.6 and 5.7, we have that all modules lie in $(**i)$ are
indecomposable modules in $\mod \Lambda_{i}$, where $1 \leq i \leq
t$. When $i\neq j$, because $S_i$ and $S_j$ do not lie in the same
block, $(**i)$ and $(**j)$ do not lie in the same block. In the
following, we will prove that all the indecomposable modules in
$\mod \Lambda _i$ lie in $(**i)$.

Let $0\neq M\in\mod\Lambda_{i}$ be indecomposable. Obviously,
$M\in\mod\Lambda$.

If $M$ is non-projective-injective with $\pd _{\Lambda}M=k$, then by
Proposition 4.8, $M\cong \Omega^{n-k}S_{j}$ for some $S_j \in
\{S_{1},S_{2},\cdots,S_{t} \}$. So $0\neq M \cong \Lambda_i M \cong
\Lambda _i \Omega^{n-k}S_{j}$ and hence $j=i$ and $M \cong \Lambda
_i\Omega^{n-k}S_{i}\cong \Omega^{n-k}S_{i}$, which lies in $(**i)$.
If $M$ is projective-injective, then we claim that $M\cong
P_{k}(S_{i})$ for some $0\leq k\leq n-1$. Since $\Lambda_{i}$ is a
block, there does not exists a projective-injective simple module in
$\mod\Lambda_{i}$. Then by Lemma 5.7 and Corollary 4.6, for any
simple module $S\in\mod\Lambda_{i}$, $S\cong \Omega^{k}S_{i}$ for
some $0\leq k\leq n$. Since $M$ is indecomposable injective, $M\cong
I^{0}(S^{'})$ for a simple module $S^{'}$ in $\mod \Lambda _i$. So
$S^{'}\cong \Omega^{k}S_{i}$ for some $1\leq k\leq n$, and hence by
Lemma 2.6 $M\cong P_{k-1}(S_{i})$ with $0\leq k-1 \leq n-1$. The
claim is proved. Thus $M$ lies in $(**i)$.

Consequently, we conclude that $\Lambda_{i}$ is a Nakayama algebra
with $\gldim\Lambda_{i}=n$ given by the quiver: $$\xymatrix{1 &
\ar[l]_{\beta _{1}} 2 & \ar[l]_{\beta _{2}} 3 & \ar[l]_{\beta _{3}}
\cdots & \ar[l]_{\beta _{n}} n+1}
$$ modulo the ideal generated by $\{
\beta_{i}\beta_{i+1}| 1\leq i\leq n-1 \}$. The proof is finished.
$\hfill\square$

\vspace{0.5cm}

\centerline{\Large \bf 6. Non-trivial maximal 1-orthogonal
subcategories}

\vspace{0.2cm}

In this section, based on [HuZ, Corollary 3.12], we will further
give a necessary condition for Auslander algebras with global
dimension 2 admitting a non-trivial maximal 1-orthogonal subcategory
in terms of the homological properties of simple modules.

\vspace{0.2cm}

{\bf Lemma 6.1} {\it Let $\Lambda$ be an Auslander algebra with
$\gldim\Lambda=2$ and $S \in \mod \Lambda$ a simple module with
$\id_{\Lambda}S=2$. Then $I^{2}(S)$ is indecomposable and
$I^{0}(S)\ncong I^{1}(S)$.}

\vspace{0.2cm}

{\it Proof.} By Lemma 3.5, we get a simple module $S^{'} \in \mod
\Lambda$ such that $\pd_{\Lambda}S^{'}=2$ and $\mathbb{D}\circ{\rm
Ext}_{\Lambda}^{2}(S^{'},\Lambda)=S$. From the minimal projective
resolution of $S^{'}$, we get an exact sequence: $$0\rightarrow
P_{0}(S^{'})^{*}\rightarrow P_{1}(S^{'})^{*}\rightarrow
P_{2}(S^{'})^{*}\rightarrow {\rm
Ext}_{\Lambda}^{2}(S^{'},\Lambda)\rightarrow 0,$$ which is a minimal
projective resolution of ${\rm Ext}_{\Lambda}^{2}(S^{'},\Lambda)$ by
[M, Proposition 4.2], where $(-)^{*}={\rm
Hom}_{\Lambda}(-,\Lambda)$. Then Applying the functor $\mathbb{D}$,
we get a minimal injective resolution of $S=\mathbb{D}\circ{\rm
Ext}_{\Lambda}^{2}(S^{'},\Lambda)$: $$0\rightarrow S\rightarrow
\mathbb{D}P_{2}(S^{'})^{*}\rightarrow
\mathbb{D}P_{1}(S^{'})^{*}\rightarrow
\mathbb{D}P_{0}(S^{'})^{*}\rightarrow 0.$$ It follows that
$I^{2}(S)\cong\mathbb{D}P_{0}(S^{'})^{*}$,
$I^{1}(S)\cong\mathbb{D}P_{1}(S^{'})^{*}$ and $
I^{0}(S)\cong\mathbb{D}P_{2}(S^{'})^{*}$. On the other hand, from
the minimal projective resolution of $S^{'}$: $$0 \rightarrow
P_{2}(S^{'}) \rightarrow P_{1}(S^{'})\rightarrow
P_{0}(S^{'})\rightarrow S^{'} \rightarrow 0,$$ we know that
$P_{0}(S^{'})$ is indecomposable and $P_{2}(S^{'})\ncong
P_{1}(S^{'})$. So our assertion follows. $\hfill{\square}$

\vspace{0.2cm}

{\bf Proposition 6.2} {\it Let $\Lambda$ be an Auslander algebra
with $\gldim\Lambda=2$. If $\Lambda$ admits a non-trivial maximal
1-orthogonal subcategory of $\mod \Lambda$, then

(1) There exists a simple module in $\mod \Lambda$ with both
projective and injective dimensions 2.

(2) There exist at least two non-injective simple modules in $\mod
\Lambda$ with projective dimension 2.}

\vspace{0.2cm}

{\it Proof.} (1) It follows from [HuZ, Corollary 3.12].

(2) By (1), there exists a non-injective simple module in $\mod
\Lambda$ with projective dimension 2. If the non-injective simple
module in $\mod \Lambda$ with projective dimension 2 is unique (say
$S$), then $\id _{\Lambda}S=2$ by (1). Since $I^{0}(S)$ and
$I^{2}(S)$ are indecomposable by Lemma 6.1, $\grade I^{0}(S)=\grade
I^{2}(S)=2$ by Lemma 3.10. Put $K=\Coker (S \hookrightarrow
I^{0}(S))$. Then $\grade K=2$ by Lemma 3.9 and so $\grade
I^{1}(S)=2$. We claim that $I^{0}(S)$ is isomorphic to a direct
summand of $I^{1}(S)$. Otherwise, since $S$ is the unique
non-injective simple module with projective dimension 2, any
non-zero indecomposable direct summand of $I^{1}(S)$ is simple by
Lemma 3.8(2). So $I^{1}(S)$ is semisimple and hence $K$ is
injective, which contradicts with $\id _{\Lambda}S=2$.

Notice that $I^{2}(S)$ is indecomposable and $\pd
_{\Lambda}I^{2}(S)=2$, so $I^{2}(S)\cong I^{0}(S)$ or $I^{2}(S)
\cong S^{'}$ for some simple module $S^{'}\in \mod \Lambda$ such
that $S \ncong S^{'}$ and $\pd _{\Lambda}S^{'}=2$. In the latter
case, we have that $\L(I^{0}(S))=\L(I^{1}(S))$. Since $I^{0}(S)$ is
isomorphic to a direct summand of $I^{1}(S)$ by the above argument,
$I^{0}(S)\cong I^{1}(S)$, which is a contradiction by Lemma 6.1.

Because $\Lambda$ is an Auslander algebra and $\pd _{\Lambda
^{op}}\mathbb{D}S=2$, it follows from Lemma 3.10 that $\grade
\mathbb{D}S=2$. Then, for any injective module $I\in \mod \Lambda$,
${\rm Ext}_{\Lambda}^{1}(I, S) \cong {\rm Ext}_{\Lambda
^{op}}^{1}(\mathbb{D}S, \mathbb{D}I)=0$. Moreover, $S
\hookrightarrow I^{0}(S)$ is left minimal, thus $K$ has no injective
direct summands by Lemma 2.5 and therefore $K$ is indecomposable by
Lemmas 6.1 and 2.1. It follows from Lemma 2.3 that $I^{1}(S)
\rightarrow I^{2}(S)$ is right minimal. So, if $I^{2}(S)\cong
I^{0}(S)$, then $I^{1}(S)$ has no simple direct summand $S^{''}$
such that $S^{''}\ncong S$ and $\pd _{\Lambda}S^{''}=2$. It yields
that $I^{1}(S) \cong [I^{0}(S)] ^{t}$ for some $t \geq 1$ and
$2\L(I^{0}(S))=t\L(I^{0}(S))+1$. It implies that $t=1$ and
$I^{0}(S)\cong I^{1}(S)$, which is a contradiction by Lemma 6.1. The
proof is finished. $\hfill{\square}$

\vspace{0.2cm}

As an immediate consequence of Proposition 6.2, we have the
following result, which gives some sufficient conditions that any
maximal 1-orthogonal subcategory of $\mod\Lambda$ (in case it
exists) is trivial for an Auslander algebra $\Lambda$.

\vspace{0.2cm}

{\bf Corollary 6.3} {\it Let $\Lambda$ be an Auslander algebra with
$\gldim \Lambda =2$. Then any maximal 1-orthogonal subcategory of
$\mod\Lambda$ (in case it exists) is trivial if one of the following
conditions are satisfied.

(1) There exists a unique simple module with projective dimension 2.

(2) There exist exactly two simple modules with projective dimension
2 and at least one of them is injective.}

\vspace{0.1cm}

From the results obtained in this paper and in [HuZ], we see that
for an $(n-1)$-Auslander algebra $\Lambda$ with $\gldim \Lambda =n$,
the properties of simple modules with projective dimension $n$ play
an important role in the study of the existence of maximal
$(n-1)$-orthogonal subcategories and the properties of $\Lambda$
admitting maximal $(n-1)$-orthogonal subcategories.

We end this section with examples to illustrate Proposition 6.2 and
Corollary 6.3.

The following example shows that there exists an Auslander algebra
$\Lambda$ with $\gldim\Lambda=2$ satisfying the condition (1) in
Proposition 6.2, but not satisfying the condition (2) in this
proposition.

\vspace{0.1cm}

{\bf Example 6.4} Let $\Lambda$ be a finite-dimensional algebra
given by the quiver:
$$\xymatrix{1 & \ar[l]_{\beta _{1}} 2 & \ar[l]_{\beta _{2}} 3
& \ar[l]_{\beta _{3}} 4 & \ar[l]_{\beta _{4}} 5}
$$
modulo the ideal generated by $\{ \beta_{1}\beta_{2},
\beta_{3}\beta_{4}\}$. Then we have

(1) $\Lambda$ is an Auslander algebra with $\gldim \Lambda =2$.

(2) All simple modules in $\mod \Lambda$ with projective dimension 2
are $S(3)$ and $S(5)$.

(3) $\id _{\Lambda}S(3)=2$ and $S(5)$ is injective.

Then by Lemma 3.7 and Proposition 6.2(2), there does not exist any
maximal 1-orthogonal subcategory of $\mod \Lambda$.

\vspace{0.1cm}

The following example shows that there exists an Auslander algebra
$\Lambda$ with $\gldim\Lambda=2$ satisfying the condition (2) in
Proposition 6.2, but not satisfying the condition (1) in this
proposition.

\vspace{0.1cm}

{\bf Example 6.5} Let $\Lambda$ be a finite-dimensional algebra
given by the quiver:
 $$\xymatrix{&6 \ar[r]^{\alpha}\ar[d]^{\gamma} & 4\ar[d]^{\beta}& \\
&5\ar[r]^{\delta}&3\ar[r]^{\lambda}\ar[d]^{\mu}&1\\
&\ \ &2&}$$ modulo the ideal generated by $\{
\beta\alpha-\delta\gamma, \mu\delta, \lambda\beta \}$. Then we have

(1) $\Lambda$ is an Auslander algebra and an almost hereditary
algebra with $\gldim \Lambda =2$.

(2) All simple modules in $\mod \Lambda$ with projective dimension 2
are $S(4)$, $S(5)$ and $S(6)$.

(3) $\id _{\Lambda}S(4)=\id _{\Lambda}S(5)=1$ and $S(6)$ is
injective.

Then by Lemma 3.7 and Proposition 6.2(1), there does not exist any
maximal 1-orthogonal subcategory of $\mod \Lambda$.

\vspace{0.2cm}

According to Examples 6.4 and 6.5, we know that the conditions (1)
and (2) in Proposition 6.2 are independent.

The following example is also related to Proposition 6.2 and
Corollary 6.3. It shows that there exists an Auslander algebra
$\Lambda$ with $\gldim\Lambda=2$ and there exists a unique simple
module $S\in \mod \Lambda$ with $\pd_{\Lambda}S=2$ and
$\id_{\Lambda}S=2$.

\vspace{0.2cm}

{\bf Example 6.6} Let $\Lambda$ be a finite-dimensional algebra
given by the quiver:
$$\xymatrix{1 \ar[r]^{\beta} & \ar@<2pt>[l]^{\alpha} 2}
$$ modulo the
ideal generated by $\beta\alpha$. Then $\Lambda$ is an Auslander
algebra with $\gldim\Lambda=2$ and the unique simple module with
projective dimension 2 is $S(2)$, and $\id_{\Lambda}S(2)=2$. Then by
Proposition 6.2(2) (or Corollary 6.3) and Lemma 3.7, there does not
exist any maximal 1-orthogonal subcategory of $\mod \Lambda$.

\vspace{0.5cm}

{\bf Acknowledgements} The research was partially supported by the
Specialized Research Fund for the Doctoral Program of Higher
Education (Grant No. 20060284002), NSFC (Grant No.10771095) and NSF
of Jiangsu Province of China (Grant No. BK2007517). The authors
thank Prof. Osamu Iyama for useful suggestions.


\vspace{0.5cm}

\begin{thebibliography}{101}

\bibitem[A]{A1}  J. L. Alperin, Local Representation Theory.
Modular Representations as an Introduction to the Local
Representation Theory of Finite Groups. Cambridge Studies in Adv.
Math. {\bf 11}, Cambridge University Press, Cambridge, 1986.

\bibitem[AsSS]{A2} I. Assem, D. Simson and A. Skowro\'nski,
Elements of the Representation Theory of Associative Algebras. Vol.
1. Techniques of Reperesentation Theory. London Math. Soc. Student
Texts {\bf 65}, Cambridge Univ. Press, Cambridge, 2006.

\bibitem[Au]{A3} M. Auslander, {\it Functors and morphisms determined
by objects.} In: Representation theory of Algebras, Proc. COnf.,
Temple Univ., Philaphia, PA,1976, in: Lecture Notes in Pure and
Appl. Math. {\bf 37}, Dekker, New York, 1978, pp.1--244.

\bibitem[AuB]{A4} M. Auslander and M. Bridger, {\it Stable module theory},
Memoirs Amer. Math. Soc. {\bf 94}, Amer. Math. Soc., Providence, RI,
1969.

\bibitem[AuR]{A5} M. Auslander and I. Reiten, {\it Applications of
contravariantly finite subcategories}. Adv. Math. {\bf 86}(1991),
111--152.

\bibitem[AuRS]{A6} M. Auslander, I. Reiten and S.O. Smol$\phi$,
Representation Theory of Artin Algebras. Corrected reprint of the
1995 original. Cambridge Studies in Adv. Math. {\bf 36}, Cambridge
Univ. Press, Cambridge, 1997.

\bibitem[B]{A7} J.E. Bj\"ork, {\it The Auslander condition on
Noetherian rings}. In: S\'{e}minaire d'Alg\`{e}bre Paul Dubreil et
Marie-Paul Malliavin, 39\`{e}me Ann\'{e}e, Paris, 1987/1988, Lect.
Notes in Math. {\bf 1404}, Springer-Verlag, Berlin, 1989,
pp.137--173.

\bibitem[EH]{A8} K. Erdmann and T. Holm, {\it Maximal n-orthogonal
modules for selfinjective algebras}. Proc. Amer. Math. Soc. {\bf
136} (2008), 3069--3078.

\bibitem[FGR]{A9} R.M. Fossum, P.A. Griffith and I. Reiten, Trivial
Extensions of Abelian Categories. Lect. Notes in Math. {\bf 456},
Springer-Verlag, Berlin, 1975.

\bibitem[GLS]{A10} C. Geiss, B. Leclerc and J. Schr\"oer, {\it rigid
modules over preprojective algebras}. Invent. Math. {\bf165} (2006),
589--632.

\bibitem[HRS]{A11} D. Happel, I. Reiten and S.O. Smal$\phi$,
{\it Tilting in Abelian categories and Quasitilted Algebras}.
Memoirs Amer. Math. Soc. {\bf 575}, Amer. Math. Soc., Providence,
RI, 1996.

\bibitem[HRi]{A12} D. Happel and C.M. Ringel, {\it Tilted algebras},
Trans. Amer. Math. Soc. {\bf 274}(1982), 399--443.

\bibitem[HuZ]{A13} Z.Y. Huang and X.J. Zhang, {\it The existence of maximal $n$-orthogonal
subcategories.}  J. Algebra (2009),
doi:10.1016/j.jalgebra.2009.01.036.

\bibitem[IS]{A14} Y. Iwanaga, H. Sato, {\it On Auslander's}
$n$-{\it Gorenstein rings}. J. Pure Appl. Algebra {\bf 106}(1996),
61--76.

\bibitem[Iy1]{A15} O. Iyama, {\it Symmetry and duality on $n$-Gorenstein
rings}. J. Algebra {\bf 269}(2003), 528--535.

\bibitem[Iy2]{A16} O. Iyama, {\it $\tau$-Categories III: Auslander
Orders and Auslander-Reiten Quivers,} Algebr. Represent. Theory {\bf
8}(2005), 601--619.

\bibitem[Iy3]{A17} O. Iyama, {\it Higher-dimensional Auslander-Reiten
theory on maximal orthogonal subcategories}. Adv. Math. {\bf
210}(2007), 22--50.

\bibitem[Iy4]{A18} O. Iyama, {\it Auslander correspondence}. Adv. Math.
{\bf 210}(2007), 51--82.

\bibitem[Iy5]{A19} O. Iyama, {\it Cluster tilting for higher Auslander algebras}.
Preprint is available at: arXiv:0809.4897 [math.RT].


\bibitem[Iy6]{A20} O. Iyama, {\it Auslander-Reiten theory revisited}.
In: Trends in Representation Theory of Algebras and Related Topics
(edited by A. Skowro\'nski), European Math. Soc. Publishing House,
Z\"urich, 2008, pp. 349--398.


\bibitem[KR]{A21} B. Keller and I. Reiten, {\it Cluster-tilted algebras are
Gorenstein and stably Calabi-Yau}. Adv. Math. {\bf 211}(2007),
123--151.

\bibitem[L]{A22} M. Lada, {\it Relative homology and maximal $l$-orthogonal modules.}
J. Algebra (2009), doi:10.1016/j.jalgebra.2009.02.015.

\bibitem[M]{A23} Y. Miyashita, {\it Tilting modules of finite
projective dimension}. Math. Z. {\bf 193}(1986), 113--146.

\end{thebibliography}
\end{document}